\documentstyle[12pt]{article}
\newlength{\abstractwidth}
\renewcommand{\thefootnote}{\fnsymbol{footnote}}
\renewcommand{\thanks}[1]{\footnote{#1}} % Use this for footnotes
\newcommand{\starttext}{ \setcounter{footnote}{0}
\renewcommand{\thefootnote}{\arabic{footnote}}}

\newcommand{\be}{\begin{equation}}
\newcommand{\bea}{\begin{eqnarray}}
\newcommand{\eea}{\end{eqnarray}} \newcommand{\ee}{\end{equation}}
 
 \def\ba{\begin{eqnarray}}
\def\ea{\end{eqnarray}}

\def\tr{{\rm tr}}

\def\log{\,{\rm log}\,}

\def\ge{\geq}
\def\le{\leq}

\def\p{\partial}

\def\[{{\bf [}}
\def\]{{\bf ]}}

\def\ddbar{i\p\bar\p}

\def\mathbb{\bf}

\newcommand{\eqref}[1]{{(\ref {#1}})}
\newcommand{\innpro}[1]{\langle {#1} \rangle}

\begin{document}
\starttext \baselineskip=18pt \setcounter{footnote}{0}
\newtheorem{theorem}{Theorem}
\newtheorem{lemma}{Lemma}
\newtheorem{corollary}{Corollary}
\newtheorem{definition}{Definition}
\newtheorem{conjecture}{Conjecture}
\newtheorem{proposition}{Proposition}

\begin{center}
{\Large \bf A NEW GRADIENT ESTIMATE FOR THE COMPLEX MONGE-AMP\`ERE EQUATION
             % A new gradient estimate for the complex Monge-Amp\`ere equation
	\footnote{Work supported in part by the National Science Foundation under grant DMS-18-55947.}}
\medskip
\centerline{Bin Guo, Duong H. Phong, and Freid Tong}

\medskip

\begin{abstract}

{\footnotesize A gradient estimate for complex Monge-Amp\`ere equations
which improves in some respects on known estimates is proved using the ABP maximum principle. }

\end{abstract}

\end{center}

\baselineskip=15pt

\section{Introduction}
\setcounter{equation}{0}
\setcounter{footnote}{0}

Gradient estimates occupy a special position in the theory of complex Monge-Amp\`ere equations. In Yau's original proof of the Calabi conjecture for compact manifolds \cite{Y}, they can be bypassed, as $C^2$ estimates can be obtained directly once $C^0$ estimates are known. But this is no longer the case for subsequent extensions of the theory. The first gradient bounds appear to be due to Hanani \cite{H}, but this paper did not seem to be widely known. More recent approaches are due to Blocki \cite{B}, P. Guan \cite{Gp}, B. Guan-Q. Li \cite{GL}, and Phong-Sturm \cite{PS}. The sharpest result to date may be \cite{PS}, which builds on the approach of \cite{B}, gives a pointwise estimate, and requires only a lower bound for the solution $\varphi$ of the equation, and not an upper bound. These features are essential for applications to boundary value problems or the case of degenerating background metrics \cite{PS, PSS}.

\smallskip
In \cite{GPT}, the authors developed a new method for establishing the classical $L^\infty$ estimates for the complex Monge-Amp\`ere equation without recourse to pluripotential theory. This method builds on works of Wang, Wang, Zhou \cite{WWZ} and particular of Chen and Cheng \cite{CC}, who introduced the idea of using an auxiliary Monge-Amp\`ere equation. The methods of \cite{GPT} turn out not just to recover the classical $L^\infty$ estimates, but to improve and widen them in many significant ways. Thus it is natural to examine their possibilities for other estimates. In this paper, we examine the case of gradient estimates. We shall show below that the methods of \cite{GPT} can recapture the sharp gradient estimate of \cite{PS, PSS}, in fact with a weaker assumption on the right hand side which may be of geometric significance.

\medskip 

Let $(X, \omega_0)$ be a compact K\"ahler manifold with or without boundary and $\varphi$ be a $\omega$-plurisubharmonic function solving the following complex Monge-Amp\`ere equation 
\begin{equation}\label{eqn:MA}
(\omega_0+\ddbar \varphi)^n = e^F \omega_0^n, 
\end{equation}
where $F\in C^\infty$. When $X$ has no boundary, we assume that $F$ satisfies the compatibility condition $\int_X e^F \omega_0^n  = \int_X \omega_0^n$. When $X$ has a smooth non-empty boundary $\p X$, we impose the boundary condition $\varphi = \phi$ on $\partial X$ for some $\phi\in C^2(\overline{X})$ with $\omega_\phi = \omega_0+ i\p\bar\p\phi$ a smooth K\"ahler metric on $\overline X$.

\smallskip

\begin{theorem}\label{thm:gradient}
Under the above conditions, we have the gradient estimate
$$|\nabla \varphi|_{\omega_0}^2 \le C e^{\lambda (\varphi - \inf_X \varphi)},$$
where $\lambda>0$ and $C$ are positive constants, depending respectively only on a lower bound for the bisectional curvature of $\omega_0$, and on $n,\omega_0, \sup_X F$, $\| \nabla F\|_{L^{2n}(e^{2F}\omega_0^n)}$, $\| \nabla_{\omega_0} \varphi\|_{L^\infty(\partial X)}$ and $\| \phi\|_{C^2(\overline X)}$.
\end{theorem}

We observe that previous results had required control of the full $L^{\infty}$ norm of the gradient of the right-hand side. With our method, we can relax this to an $L^{2n}$-control. We illustrate later an application of this improvement. 

\section{Proof of the Theorem}
\setcounter{equation}{0}

 By replacing $\varphi$ and $\phi$ respectively by
$\varphi - \inf_X \varphi$ and $\phi - \inf_X \phi$, we may assume $\inf_X \varphi = 0$. Let $\omega = \omega_0+ \ddbar \varphi$ be the K\"ahler metric associated with the complex Monge-Amp\`ere  equation \eqref{eqn:MA}. 

\begin{lemma}\label{lemma gradient}
The following equation holds
\bea\label{eqn:gradient}
\Delta_\omega |\nabla \varphi|^2_{\omega_0} = 2 Re\innpro{\nabla F, \bar\nabla \varphi}_{\omega_0} + g^{i\bar j} g_{0}^{k\bar l} ( \varphi_{ki}\varphi_{\bar j \bar l} + \varphi_{k\bar j} \varphi_{i\bar l}  ) + g^{i\bar j} R(g_0)_{i\bar j k\bar l} \varphi_p \varphi_{\bar q} g^{k\bar q}_0 g_0^{p\bar l}
\eea
where $\omega = (g_{i\bar j})$, $\omega_0 = ((g_0)_{i\bar j})$, $\varphi_{k i } = (\nabla_{\omega_0} \nabla_{\omega_0} \varphi)_{ki}$ are the second covariant derivatives with respect to $\omega_0$, and $R(g_0)_{i\bar j k\bar l}$ is the bisectional curvature of $\omega_0$.
\end{lemma}
The proof of Lemma \ref{lemma gradient} is a standard calculation, so we omit the details. Let $-K$ be a lower bound of the bisectional curvature $R(g_0)$.
From the equation \eqref{eqn:gradient} we have
\bea\nonumber
\Delta_\omega |\nabla \varphi|^2_{\omega_0} \ge  2 Re\innpro{\nabla F, \bar\nabla \varphi}_{\omega_0} + g^{i\bar j} g_{0}^{k\bar l} ( \varphi_{ki}\varphi_{\bar j \bar l} + \varphi_{k\bar j} \varphi_{i\bar l}  )  - 2 K \tr_{\omega} \omega_0 |\nabla \varphi|^2_{\omega_0}
\eea
Denote $H = e^{-\lambda \varphi} |\nabla \varphi|_{\omega_0}^2$ for $\lambda = 2 K + 10$. We calculate at an arbitrary point $x\in X$, and choose a normal coordinates system for $\omega_0$ such that $\omega$ is diagonal at $x$. 
\bea
\Delta_\omega H& =& \Delta (e^{-\lambda \varphi} |\nabla \varphi|^2_{\omega _0}) \nonumber \\
 \nonumber & = &e^{-\lambda \varphi} \Delta |\nabla \varphi|_{\omega_0}^2 + |\nabla \varphi|^2_{\omega_0} \Delta(e^{-\lambda\varphi }) - 2 \lambda e^{-\lambda\varphi}Re \innpro{ \nabla \varphi, \bar \nabla |\nabla \varphi|^2_{\omega_0}  }_\omega\\
&\ge & \nonumber e^{-\lambda \varphi} \Big(2 Re\innpro{\nabla F, \bar\nabla \varphi}_{\omega_0} + g^{i\bar j} g_{0}^{k\bar l} ( \varphi_{ki}\varphi_{\bar j \bar l} + \varphi_{k\bar j} \varphi_{i\bar l}  )  - 2 K \tr_{\omega} \omega_0 |\nabla \varphi|^2_{\omega_0} \Big)\\
&\quad &\nonumber+ |\nabla \varphi|_{\omega_0}^2 e^{-\lambda \varphi} \Big( - \lambda n + \lambda \tr_\omega \omega_0 + \lambda^2 |\nabla\varphi|_\omega^2   \Big) - 2 \lambda e^{-\lambda\varphi}Re \innpro{ \nabla \varphi, \bar \nabla |\nabla \varphi|^2_{\omega_0}  }_\omega
\eea
The last term on the right hand side is
\bea
 \nonumber && - 2 \lambda e^{-\lambda\varphi}Re \innpro{ \nabla \varphi, \bar \nabla |\nabla \varphi|^2_{\omega_0}  }_\omega \\ &\nonumber = & - 2\lambda e^{-\lambda\varphi } Re\big( g^{i\bar i } \varphi_i (\varphi_k \varphi_{\bar k})_{\bar i} \big)\\
 &\nonumber = & -2\lambda e^{-\lambda\varphi} Re \big( g^{i\bar i} \varphi_i \varphi_{k\bar i} \varphi_{\bar k} + g^{i\bar i}  \varphi_i \varphi_k \varphi_{\bar k \bar i}     \big)\\
 &\nonumber\ge & -2 \lambda e^{-\lambda \varphi} g^{i\bar i} \varphi_i \varphi_{k\bar i} \varphi_{\bar k} { - \lambda^2 e^{-\lambda\varphi} g^{i\bar i} \varphi_i \varphi_{\bar i} |\nabla \varphi|^2_{\omega_0}-  e^{-\lambda \varphi} g^{i\bar i} \varphi_{k i }\varphi_{\bar k \bar i}}\\
 & \nonumber= & -2\lambda e^{-\lambda\varphi } g^{i\bar i} \varphi_i \varphi_{\bar i} ( g_{i\bar i} - 1  ) { - \lambda^2 e^{-\lambda\varphi} g^{i\bar i} \varphi_i \varphi_{\bar i} |\nabla \varphi|^2_{\omega_0}-  e^{-\lambda \varphi} g^{i\bar i} \varphi_{k i }\varphi_{\bar k \bar i}} \\
 & = & \nonumber - 2\lambda e^{-\lambda \varphi} |\nabla \varphi|_{\omega_0}^2 + 2\lambda e^{-\lambda\varphi} |\nabla \varphi|^2_{\omega} { - \lambda^2 e^{-\lambda\varphi} |\nabla \varphi|^2_\omega |\nabla \varphi|^2_{\omega_0}-  e^{-\lambda \varphi} g^{i\bar i} \varphi_{k i }\varphi_{\bar k \bar i}}.
\eea
where we applied the Cauchy-Schwarz inequality.
We thus obtain  
\bea
\nonumber
\Delta_\omega H& \ge 2 e^{-\lambda \varphi} Re\innpro{\nabla F, \bar \nabla \varphi}_{\omega_0} + (\lambda - 2K) H\tr_\omega \omega_0  - \lambda(n+2) H.
\eea
Note that this inequality holds at any point of $X$, since it is independent of the choice of normal coordinates. Let $\alpha>1$ be a positive constant. We calculate
\bea\nonumber
\Delta_\omega H^\alpha& = &\alpha H^{\alpha -1} \Delta H + \alpha (\alpha - 1) H^{\alpha - 2} |\nabla H|^2_{\omega}\\
&\ge & \nonumber \alpha H^{\alpha - 1} \Big( e^{-\lambda \varphi} Re\innpro{\nabla F, \bar \nabla \varphi}_{\omega_0}  
+(\lambda - 2K) H\tr_\omega \omega_0  - \lambda(n+2) H  \Big)
\nonumber\\
&& + \alpha (\alpha -1) H^{\alpha -2} |\nabla H|_\omega^2.
\eea
Since $X$ is compact, we can assume $H$ attains its maximum at a point $x_0$, with $H(x_0) =: M>0$. We may suppose $x_0$ lies in the interior of $X$, otherwise we are done. Since $\omega_0$ is smooth up to $\partial X$, we may assume $(\overline X,\omega_0)$ isometrically embeds to another K\"ahler manifold $(\hat X, \hat \omega_0) $ as a compact subset\footnote{We can alternatively cover $\partial X$ by finitely many Euclidean half balls, and apply similar calculations.}.% 
 Let $r>0$ be the injectivity radius of the Riemannian manifold  $(\hat X,\hat \omega_0)$. Without loss of generality we may identify the metric ball $B_{g_0}(x_0, r)$ with an open domain in the Euclidean space $\mathbb C^n$, where we denote $B_{g_0}(x_0, r) = \{x\in X | d_{g_0}(x,x_0)<r\}$. We will apply a trick of Chen-Cheng \cite{CC}. Let $\theta =  \min\{\frac{1}{10n C_0}, \frac{r^2}{10n C_0}  \}$ be a given constant (where $C_0>1$ depends only on $\omega_0$) and choose an auxiliary function $\eta$ such that $\eta = 1$ on $B_{g_0}(x_0, r/2)$ and $ \eta = 1-\theta$ on $B_{g_0}(x_0,r)\backslash B_{g_0}(x_0, 3r/4)$, and $\eta \in [1-\theta, 1]$ in the annulus between. We also have (this $\eta$ may be chosen as $\hat\eta ( \frac{d_0(x)^2}{r^2}  )$ where $d_0$ is a smoothing of the $g_0$-distance to $x_0$ and $\hat \eta$ is some appropriate function on $\mathbb R$)
$$|\nabla \eta|_{g_0}^2 \le \frac{C_0\theta^2}{r^2},\quad |\nabla^2 \eta|_{g_0}\le \frac{C_0 \theta}{r^2}$$
We calculate as follows
\bea\nonumber
\Delta_\omega (\eta H^\alpha) = \eta \Delta H^\alpha + 2\alpha H^{\alpha - 1} Re\innpro{\nabla \eta,  \bar \nabla H }_\omega + H^\alpha \Delta_\omega \eta.
\eea
Note that the last term satisfies $$H^\alpha \Delta _\omega \eta \ge - C_0 \frac{\theta}{r^2} H^\alpha  \tr_\omega \omega_0,$$
and the middle term is
\bea
 2\alpha H^{\alpha - 1} Re\innpro{\nabla \eta,  \bar \nabla H }_\omega  & \ge & \nonumber - 2 \alpha H^{\alpha - 1} |\nabla H|_\omega |\nabla \eta|_\omega \\
 & \ge & \nonumber - \frac{\alpha (\alpha -1)}{2} H^{\alpha - 2} |\nabla H|_\omega^2 - \frac{2\alpha }{\alpha - 1} H^\alpha |\nabla \eta|_\omega^2\\
  & \ge & \nonumber \underbrace{- \frac{\alpha (\alpha -1)}{2} H^{\alpha - 2} |\nabla H|_\omega^2}_{\mbox{controlled by the last term in $\eta \Delta H^\alpha$}} - \frac{2\alpha }{\alpha - 1} H^\alpha \frac{C_0\theta^2}{r^2} \tr_\omega \omega_0
\eea
Combining the above inequalities we get
\bea\nonumber
\Delta (\eta H^\alpha)&\ge& \alpha \eta H^{\alpha - 1} e^{-\lambda \varphi} \innpro{\nabla F, \bar \nabla \varphi}_{\omega_0} + 
 \Big(\alpha \eta(\lambda - 2K) - \frac{C_0 \theta}{r^2} - \frac{2\alpha}{\alpha - 1} \frac{C_0 \theta^2}{r^2} \Big) H^{\alpha} \tr_\omega \omega_0
 \nonumber\\
 && - \lambda \alpha (n+2) H^\alpha.
\eea
Note that $\eta\ge 9/10$. We can choose $\alpha =2$. Together with the choice of $\theta$ and $\lambda$, the middle term of the right hand side of the above inequality is nonnegative, so
\bea
\Delta_\omega (\eta H^\alpha) \ge -\alpha \eta H^{\alpha - \frac 1 2 } e^{-\lambda \varphi/2} |\nabla F|_{\omega_0} - \lambda \alpha( n+2) H^\alpha
\eea

We may assume $(1-\theta) M^\alpha \ge \sup_{\partial X} |\nabla \varphi|^\alpha_{\omega_0}$, otherwise we are done. 
Applying the ABP maximum principle to the function $\eta H^\alpha$ on the ball $B_{g_0} (x_0, r)$, we obtain (with $B_0 = B_{g_0}(x_0,r)$)
\bea
M^\alpha & = &\nonumber \sup_{B_0} (\eta H^\alpha)\\
&\le&\nonumber \sup_{\partial B_0}( \eta H^\alpha )+ C(n,\omega_0) r \Big( \int_{B_0}  \frac{  \big[\alpha \eta H^{\alpha - \frac 1 2 } e^{-\lambda \varphi/2} |\nabla F|_{\omega_0}  + \lambda \alpha( n+2) H^\alpha \big]^{2n}   }{e^{-2 F}} \omega_0^n\Big)^{1/2n}\\
&\le&\nonumber \sup_{\partial B_0} \eta H^\alpha + C(n,\omega_0) r\Big[ \Big( \int_{B_0} H^{2n \alpha} \omega_0^n \Big)^{1/2n} + M^{\alpha - 1/2} \big(\int_{B_0} e^{2F} |\nabla F|_{\omega_0}^{2n} \omega_0^n \big)^{1/2n}   \Big]\\
&\le&\nonumber \sup_{\partial B_0} \eta H^\alpha + C(n, \omega_0) r\Big[ M^{\alpha (1-\frac 1{2n})} \Big( \int_{B_0} H\omega_0^n  \Big)^{1/2n} + M^{\alpha - 1/2}  \big(\int_{B_0} e^{2F} |\nabla F|_{\omega_0}^{2n} \omega_0^n \big)^{1/2n}  \Big]\\
&\le&\nonumber (1-\theta) M^\alpha+ C(n,\omega_0, F) r\Big[ M^{\alpha (1-\frac 1{2n})} + M^{\alpha - 1/2}    \Big],
\eea 
where the last constant $C(n,\omega_0, F)$ depends on $\| |\nabla F|_{\omega_0} \|_{L^{2n}(X, e^{2F}\omega_0^n)}$ and $\sup_X F$. Note that $\theta  > c_0>0$ for some constant $c_0$ depending only on $\omega_0$.  We conclude that 
$$c_0 M^\alpha \le \theta M^\alpha\le C(n,\omega_0, F) r ( M^{\alpha - 1/2} + M^{\alpha (1-\frac 1{2n})}   ),$$
from which we derive $M\le C(n,  \omega_0,F)$, since the RHS are powers of $M$ with degree smaller than $\alpha$. Finally in the estimate above we implicitly use the uniform bound on $\int_X H \omega_0^n$, which follows from the lemma below. It is because of this lemma that we need the $C^2$-bound of the boundary value $\phi$.

\begin{lemma}
	We have $\int_X H \omega_0^n = \int_X e^{-\lambda \varphi} |\nabla \varphi|_{\omega_0}^2 \omega_0^n \le C(n, \omega_0,\omega_\phi)$.
	
\end{lemma}
{\em Proof.} From the equation $\omega^n = e^F \omega_0^n$, we obtain
$$e^F  \omega_0^n - \omega_\phi^n = \omega^n - \omega_\phi^n = \ddbar( \varphi - \phi) \wedge (\omega^{n-1}+ \cdots + \omega_\phi^{n-1}  ),$$
Multiplying both sides by $e^{-\lambda \varphi + \lambda \phi} - 1$ and applying integration by parts, we can write the right hand side as
\bea\nonumber
&& \int_X (e^{-\lambda \varphi + \lambda \phi} - 1) \ddbar (\varphi-\phi) \wedge (\omega^{n-1} + \cdots + \omega_0^{n-1})\\
& = & \int_X \lambda e^{-\lambda (\varphi - \phi)} \partial (\varphi-\phi) \wedge \bar \partial (\varphi - \phi) \wedge (\omega^{n-1} + \cdots + \omega_\phi^{n-1}) \nonumber\\
&\ge&\nonumber   \int_X \lambda e^{-\lambda (\varphi - \phi)} \partial (\varphi-\phi) \wedge \bar \partial (\varphi-\phi) \wedge \omega_\phi^{n-1} \\
&\ge& \nonumber C \int_X e^{-\lambda\varphi} |\nabla \varphi|_{\omega_0}^2 \omega_0^n - C(\phi,\omega_0),
\eea
since $\omega_\phi$ is equivalent to $\omega_0$ by assumption. 
On the other hand the left hand side can be bounded as follows,
$$\int_X (e^{-\lambda (\varphi - \phi)} - 1) (e^F\omega_0^n - \omega_\phi^n) \le \int_X e^{\lambda \phi} (e^F \omega_0^n+ \omega_\phi^n) \le C(n,\omega_0, \phi,\sup_X F).$$
This proves the lemma, and 
the proof of the theorem is complete. 

\section{Application}
\setcounter{equation}{0}

We can give now an application of the improved gradient estimates: 

\begin{corollary}
	Let $X$ be a compact K\"ahler manifold, $f\geq 0$ be a smooth function with $\int_X f\omega_0^n = \int_X \omega_0^n$, and assume that $f$ satisfies 
	\begin{equation}
		\int_X \frac{|\nabla f|^{2n}}{f^{2n-2}}\omega_0^n<\infty
	\end{equation}
	Then the solution of the complex Monge-Amp\`ere equation 
	\begin{equation}
		(\omega_0+\ddbar \varphi)^n = f\omega_0^n
	\end{equation} 
	is Lipschitz continuous. 
\end{corollary}
{\em Proof.} Let $f_{k}$ be a regularization of $f$ chosen with the following properties: $\int_X f_{k}\omega_0^n = \int_X f\omega_0^n$, $f_{k}>0$ for $k>0$, $f_k$ converges to $f$ smoothly as $k\to \infty$, and moreover in a small neighborhood of the vanishing locus of $f$, we require that $f_k = f + \frac{1}{k}$. It's not hard to see that such a regularization can be arranged. 

Then by our choice of $f_k$, we observe that in a neighborhood of the vanishing locus of $f$, we have
\begin{equation}
	\frac{|\nabla f_k|^{2n}}{f_k^{2n-2}} = \frac{|\nabla f|^{2n}}{f_k^{2n-2}} \leq \frac{|\nabla f|^{2n}}{f^{2n-2}} 
\end{equation}
and hence $\int_X \frac{|\nabla f_k|^{2n}}{f_k^{2n-2}}\omega_0^n$ is uniformly bounded. By Theorem~\ref{thm:gradient}, we know that the sequence of solutions of the Monge-Amp\`ere equations
\begin{equation}
	(\omega_0+\ddbar \varphi_k)^n = f_k\omega^n
\end{equation}
has uniform $C^1$ bounds independent of $k$. By the stability of complex Monge-Amp\`ere equations \cite{K, GPT1}, the solutions $\varphi_k$ converge uniformly to $\varphi$, hence $\varphi$ must be Lipschitz. 

\bigskip 
Our theorem applies for example, when $f$ has isolated zeroes, near which $f$ is asymptotically $f(z)\sim \frac{1}{\big|\log |z|^2\big|}$ or $f(z) \sim |z|^{\varepsilon_0}$ for $\varepsilon_0>0$. Observe that $\|\nabla f^{1/n}\|_{L^\infty}$ is not finite, so the usual gradient estimate in \cite{B,PS} do not apply. However we still have $\int_X |f|^2 |\nabla \log f|^{2n} \omega_0^n <\infty$, hence our result applies.

\bigskip

\noindent Department of Mathematics \& Computer Science, Rutgers University, Newark, NJ 07102 USA

\noindent bguo@rutgers.edu

\bigskip

\noindent Department of Mathematics, Columbia University, New York, NY 10027 USA

\noindent phong@math.columbia.edu, tong@math.columbia.edu

\end{document}